# Accurate Evaluation of Polynomials

## Brian M. Sutin

## Claremont Graduate University
School of Mathematical Sciences
710 N. College Ave.
Claremont, CA 91711

## 0. Abstract

For a large class of polynomials, the standard method of polynomial evaluation, Horner's method, can be very inaccurate. The alternative method given here is on average 100 to 1000 times more accurate than Horner's Method. The number of floating point operations is twice that of Horner's method for a single evaluation. For repeated evaluations at nearby points, the number of floating point operations is only doubled for the first evaluation, and is the same as Horner's Method for all following evaluations. This new method is tested with random polynomials.

## 1. Introduction

Quickly and accurately finding zeros of polynomials is important for computer graphics (CG) and computer aided design (CAD). However, many polynomials cannot be evaluated accurately enough near a zero to give a sufficiently accurate result. A recent survey [2] of methods for evaluation of polynomials compares 3 different methods of improving accuracy to Horner's method. The fastest of these alternate methods is 4 to 20 times slower that Horner's method over the range of polynomials considered. The method presented here is in practice the same speed as Horner's method for zero finding while increasing the accuracy by approximately a factor of up to 1,000. This is done by using a nearby polynomial that can be more accurately evaluated.

## 2. Floating Point

Modern computers almost all represent real numbers as floating point, as specified in the IEEE 754 standard [3]. The IEEE representation expresses each number as a mantissa between ½ and 1, truncated to some number of binary bits, multiplied by a binary exponent. For a signed mantissa M and an exponent E, the resulting expressed floating point number is $M\ 2^E$.

Floating point operations have two causes of error that apply to Horner's method. The first is multiplication; when two numbers A and B are multiplied, $M_A\ M_B = M_{AB}$ should have approximately twice as many bits, but is truncated to fit within the standard. The second cause of error is addition; if the exponents are different, then bits from mantissa of the number with the smaller exponent will be discarded.

3. Horner's Method

Horner's method is an algorithm for evaluating single-variable polynomials by recursion. The algorithm can be written as follows:

```
# Horner's method to evaluate a polynomial at a point
# Inputs are the polynomial coefficients P_{0...n} and a point x.
  H = 0
  for i from n to 0
        H = H x + P_i
# Output is H, the value of the polynomial evaluated at x.
```

Horner's method has a multiplication and addition at every step, along with the associated floating point error. These errors are particularly troublesome for zero-finding algorithms. When the polynomial is evaluated near a zero, some catastrophic cancellation (subtracting two nearby numbers) is guaranteed to occur.

## 4. Other Methods of Accurately Evaluating Polynomials

Hoffmann *et al* [2] survey various methods for accurate evaluation of polynomials. Three methods are compared. The first, *accumulation*, uses an extended precision accumulator to more accurately compute a sum of products, or inner product. The second method uses *distillation*, where a sum of floating point numbers is accurately represented as equivalent to a second set of disjoint floating point numbers, where disjoint here means that the sum of any two floating point numbers always results in the larger magnitude number after floating point truncation. The third method is distillation combined with an iterative improvement step. The distillation methods are both considerably slower than accumulation.

For the example polynomials used in Hoffmann *et al*, accumulation was usually slower than Horner's method by about a factor of 4 to 5, but sometimes by as much as a factor of more than 20. In contrast, the method given here is a factor of two slower if the polynomial is to be evaluated only once. For multiple nearby evaluates, such as for zero finding, the extra overhead becomes negligible. On the other hand these other methods may be more accurate for some extremely difficult-to-evaluate polynomials, where the method proposed here may be as poor as Horner's method.

## 5. Accurate Evaluation Using a Nearby Polynomial

We want to accurately evaluate the polynomial

$$P(x) = \sum_{i=0}^{n} P_i x^i \quad . \tag{1}$$

First we find a nearby polynomial that we can evaluate exactly, even with finite decimal place arithmetic. By exact evaluation, we mean if the polynomial is evaluated using finite-precision floating point and Horner's Method, then no precision will be lost relative to a computation that is exact. Our exact polynomial, with coefficients to be defined later, is

$$\hat{P}(\hat{x}) = \sum_{i=0}^{n} \hat{P}_i \hat{x}^i ,\qquad(2)$$

where $\hat{x}$ is close to $x$ and $\hat{P}$ is close to $P$. Now rearrange $P(x)$,

$$P(x) = P(\hat{x}) + \left( P(x) - P(\hat{x}) \right) .\qquad(3)$$

The first term is a constant. Then

$$P(x) = P(\hat{x}) + \sum_{i=1}^{n} P_i (x^i - \hat{x}^i) .\qquad(4)$$

Each of the terms on the right is divisible by $x - \hat{x}$, so this expression can be divided out. This is similar to the trick used in Enenkel and Keras (2004).

$$P(x) = P(\hat{x}) + (x - \hat{x}) \sum_{i=1}^{n} P_i \sum_{j=0}^{i-1} x^j \hat{x}^{i-1-j} .\qquad(5)$$

Now swapping the summations gives a polynomial in $x$,

$$P(x) = P(\hat{x}) + (x - \hat{x}) \sum_{j=0}^{n-1} \left[ \sum_{i=j+1}^{n} P_i \hat{x}^{i-1-j} \right] x^j .\qquad(6)$$

Defining the constants of this polynomial as $\{C_j\}$,

$$C_j = \sum_{i=j+1}^{n} P_i \hat{x}^{i-1-j} = \sum_{i=0}^{n-1-j} P_{j+1+i} \hat{x}^i ,\qquad(7)$$

the polynomial becomes

$$P(x) = C_{-1} + (x - \hat{x}) \sum_{j=0}^{n-1} C_j x^j .\qquad(8)$$

In order to increase accuracy, the computation for $\{C_j\}$ can be broken up into two parts,

$$C_j = \sum_{i=0}^{n-1-j} \hat{P}_{j+1+i} \hat{x}^i + \sum_{i=0}^{n-1-j} (P_{j+1+i} - \hat{P}_{j+1+i}) \hat{x}^i .\qquad(9)$$

There are two things to note here. The first summation is computed exactly, since the partial sums are just the partial sums from Horner's Method. By assumption, this computation is exact. The second summation is composed of terms that are small, since $\hat{P}$ is assumed to be close to P.

## 6. Finding a Nearby Polynomial

Each step of Horner's Method may be written out as

$$H_n = H_{n-1} x + P_n \ . \tag{10}$$

The loss in precision takes place a two places: the multiplications and the additions. The loss at the multiplications can be eliminated by truncating the lower bits of x and assuring that $H_{n-1}$, the partial sum from the previous iteration, is computed to be a truncated value. For the addition, the actual number of mantissa bits lost is

$$B_{lost} = \left| \log_2 |P_n| - \log_2 |H_{n-1} x| \right| \ . \tag{11}$$

Given $P_n$ and $H_{n-1}x$, $B_{lost}$ can be computed to sufficient accuracy with no floating point operations by subtracting the floating point exponents.

Consider a polynomial that serendipitously has all the $B_{lost}$ near zero. If the floating point mantissa has M bits, then define a function $T_{M/2}(x)$ that truncates the argument to M/2 bits. An algorithm to find a nearby $\hat{P}$ and $\hat{x}$ might be

```
# Algorithm to compute a nearby polynomial for a class of easy polynomials
# Inputs are the polynomial coefficients P_{0...n} and a point x in a region of interest.
    x̂ = T_{M/2}(x)
    H = 0
    for i from n to 0
        S = T_{M/2}(H x̂)
        H = T_{M/2}(S + P_i)
        P̂_i = H - S
# Outputs are a set of polynomial coefficients P̂_{0...n} near P and a point x̂ near x.
```

The resulting polynomial can be exactly evaluated at $\hat{x}$. At every multiplication step, H and $\hat{x}$ each have M/2 bits, so the product can be computed exactly. The $\hat{P}$ is chosen to exactly cancel out the lower M/2 bits of the product, so that the result of the addition, H, is once again only M/2 bits.

Note that the split of accuracy between M/2 bits for the nearby polynomial and M/2 bits for the nearby x was completely arbitrary. Instead we could have chosen 3M/4 and M/4 or *vice versa*. The accuracy of any given split could be more or less accurate, depending on the polynomial. In practice, numerical experiments show that M/2 appears to be optimal for random polynomials.

## 7. Problems Computing the Approximating Polynomial

If the exponents for the arguments of an addition are sufficiently far apart, then the above algorithm will fail to work as designed. For example, the polynomial

$$P(x) = 10^{-30} x^2 + 10^{30} x - 10^{30} \tag{12}$$

evaluated near x = 1 requires a machine floating point precision of at least 200 bits, otherwise the quadratic term will have no effect on the result. For polynomials like this, any $\hat{P}$ either won't be nearby $P$ or $\hat{P}$ cannot be evaluated exactly. In practice, numerical experiments with random polynomials show that choosing $\hat{P}$ nearby $P$ gives better results, as opposed to choosing a far away $\hat{P}$ that can be evaluated exactly.

Define the function $B(z)$, returning the actual number of bits in the mantissa of a floating point number. Since all numbers except zero have at least one bit, the first bit is suppressed in IEEE 754. For example, $B(2^n) = 0$, while $B(3) = 1$, and $B(14) = 2$. Next define the function $E(z)$, returning the exponent of a floating point number. Since the mantissa is defined to be in [0.5,1), some examples are $E(2^n) = n\text{-}1$, $E(3) = 2$, and $E(14) = 4$. Let M be the maximum number of bits possible in the floating point mantissa. Under IEEE 754 single-precision floating point, M is 23. The gap in exponents for the $n^{th}$ addition that causes an issue with computing $\hat{P}$ is

$$\Delta_n = E(P_n) - E(H_{n+1}) - E(\hat{x}) \ . \tag{13}$$

The conditions on $\Delta_n$ that are required to have an exact $\hat{P}$ are given in the following table.

|  | $\Delta_n > 0$ | $\Delta_n < 0$ |  |
|---|---|---|---|
| Condition on Input | $B(H_{n+1}) \leq M - B(\hat{x}) - \Delta_n$ | $B(H_{n+1}) \leq M - B(\hat{x})$ | (14) |
| Condition on Output | $B(H_n) \geq \Delta_n$ | $B(H_n) \geq -\Delta_n$ |  |

For a relatively benign polynomial, all of these conditions are simple enough to meet using $T_{M/2}()$. In the case of more extreme polynomials, finding a truncation function and $\hat{P}$ that meets the above conditions at every step can be extremely time consuming or impossible. Instead, the simple heuristic algorithm presented below is used to find a $\hat{P}$ that seems to work better than any small nearby variations of the fixed parameters.

> \# Compute the number of bits to truncate given a difficult polynomial
> \# Inputs are the polynomial coefficients $P_{0...n}$ and a point $x$ in a region of interest.
> $\hat{x} = T_{M/2}(x)$
> $\hat{M} = max\left(M - 1 - \underset{n}{\operatorname{argmax}} \Delta_n, (M+1)/3\right)$
> \# Output is a reduced number of mantissa bits, $\hat{M}$ .

This algorithm takes the same amount of time as Horner's method, so for some difficult polynomials, the time required to evaluate the polynomial the first time near a point is three times Horner's method: once to compute $\hat{M}$, once to compute the nearby polynomial, and then once for each evaluation.

The result can again be improved somewhat in the case where $x$ has bits with value zero near the point of truncation. Shifting these unused bits from $\hat{x}$ to $\hat{P}$ improves the accuracy. For example, for a polynomial being evaluated at $x = 1$, $\hat{P}$ could use the entire $\hat{M}$ bits. This looks like

> \# Improve truncation values for serendipitous values of $x$
> \# Inputs are $x$ and $\hat{M}$
> $R = \hat{M}/2$

```
while( T_{M-R}(x) = T_{M-R-1}(x) )
    R = R + 1
    # Output is a reduced number of mantissa bits, M̂, and the number of truncation bits, R.
```

Next $x$ is truncated to $\hat{x} = T_{M-R}(x)$, and the function $T_R()$ is used in the algorithm of section 5. This extra step can be done using no floating point operations, and so requires negligible extra run time.

## 8. Testing with Random Polynomials

Jenkins & Traub (1975) give three methods for generating random polynomials likely to stretch the abilities of a zero finding algorithm. The one chosen for this paper is from section 7 (v). The mantissa and exponent of each root are drawn from independent uniform distributions. The final polynomial is then computed by multiplying together linear terms based on each root using floating point arithmetic. This method has the simplification that the roots of the final polynomial are crudely known in advance. The uniform distribution for the mantissa is U{ (0.5,1) + (-1,-0.5) }. Let the final polynomial be of order N, the maximum floating point exponent be F, and the 'difficulty parameter' be D. Then the exponent is drawn from the uniform distribution U(-F/ND,F/ND). Under IEEE 754 single precision, F is 127. If D is 1, then the polynomials are as extreme as possible without overflowing or underflowing any of the terms in the final polynomial. A large D makes more innocuous polynomials. In the limit as D goes to infinity, the roots all end up uniformly distributed in ±(0.5,1).

Figures 1 and 2 show some simulations in single precision floating point arithmetic. The heuristic improvements from section 7 were all used. Figure 1 shows $8^{th}$ order polynomials with varying difficulty parameters. 128 polynomials were evaluated for a total of 1024 roots. Figure 2 shows polynomials of various orders, all with the highest difficulty level of D = 1. All the roots of each polynomial were evaluated, with the number of polynomials chosen to give 1024 points on each plot.

Directly comparing the accuracy of methods for evaluating near the zero is difficult because any direct ratio is likely to overflow. Instead here both methods are normalized by dividing each by the worst-case expected numerical accuracy, as if assuming that the coefficients of some original polynomial had all been rounded to single precision. Explicitly, this value is

$$E_{max} = \left| \frac{dP(x)}{dx} \right| \frac{eps}{2} + \sum_n \left| n P_n x^{n-1} \right| \frac{eps}{2} \quad . \tag{15}$$

Here *eps* is is the machine precision, the largest floating point number such that, when added to 1, results again in 1.

## 9. Conclusions

An algorithm has been presented that evaluates polynomials more accurately that Horner's method. The runtime is 1/(1+1/N) compared to Horner's method for N evaluations. Although the algorithm as presented is not guaranteed to have better accuracy than Horner's method, the typical accuracy improvement is a factor of 100 to 1000 for single precision floating point arithmetic. There are, no doubt, better ways to find a nearby polynomial and thus improve the accuracy of this algorithm, but the method given here is both simple and fast. This algorithm nicely fills a niche between Horner's method and much more complicated and slower algorithms.

I would like to thank Ali Nadim for many useful suggestions. This work has been entirely supported by the United Technology Corporation (UTC) Employee Scholarship Program.

## Appendix - Deflation

For finding all of the zeros of a polynomial, the next step after find a zero is often to deflate the polynomial, or in other words, divide out the known root so that the next root may be found. The polynomial can be deflated in terms of $\{C_j\}$. The deflated polynomial $S$ is defined by the equation

$$P(x) = P(r) + (x - r)S(x) \quad . \tag{16}$$

where $r$ is the recently computed root. Since the iterative procedure to find $r$ has some stopping criterion, $P(r)$ is very small but most likely not exactly zero. Now trivially

$$P(x) = P(r) + (P(x) - P(r)) \quad . \tag{17}$$

Now using equation (5.8), $P(x)$ and $P(r)$ on the right hand side can be rewritten in terms of $\{C_j\}$ as well:

$$P(x) = P(r) + (x - \hat{x})\sum_{j=0}^{n-1} C_j x^j - (r - \hat{x})\sum_{k=0}^{n-1} C_k r^k \quad . \tag{18}$$

Breaking up $(x - \hat{x})$ as $(x - r) + (r - \hat{x})$ gives

$$P(x) = P(r) + (x - r)\sum_{j=0}^{n-1} C_j x^j + (r - \hat{x})\sum_{k=1}^{n-1} C_k (x^k - r^k) \quad . \tag{19}$$

The rightmost term is divisible by $x - r$, so this divisor can be pulled out

$$P(x) = P(r) + (x - r)\sum_{j=0}^{n-1} C_j x^j + (r - \hat{x})\sum_{k=1}^{n-1} C_k (x - r)\sum_{i=0}^{k-1} x^i r^{k-1-i} \tag{20}$$

and by swapping the summations,

$$P(x) = P(r) + (x - r)\sum_{j=0}^{n-1} C_j x^j + (x - r)\sum_{i=0}^{n-2} x^i (r - \hat{x}) \sum_{k=i+1}^{n-1} C_k r^{k-1-i} \quad . \tag{21}$$

Equating this to the definition for the deflated polynomial gives the constants of the deflated polynomial S(x) as

$$S_{j<n-1} = C_j + (r - \hat{x}) \sum_{k=j+1}^{n-1} C_k r^{k-1-j} \tag{22}$$

$$S_{n-1} = C_{n-1} \quad . \tag{23}$$

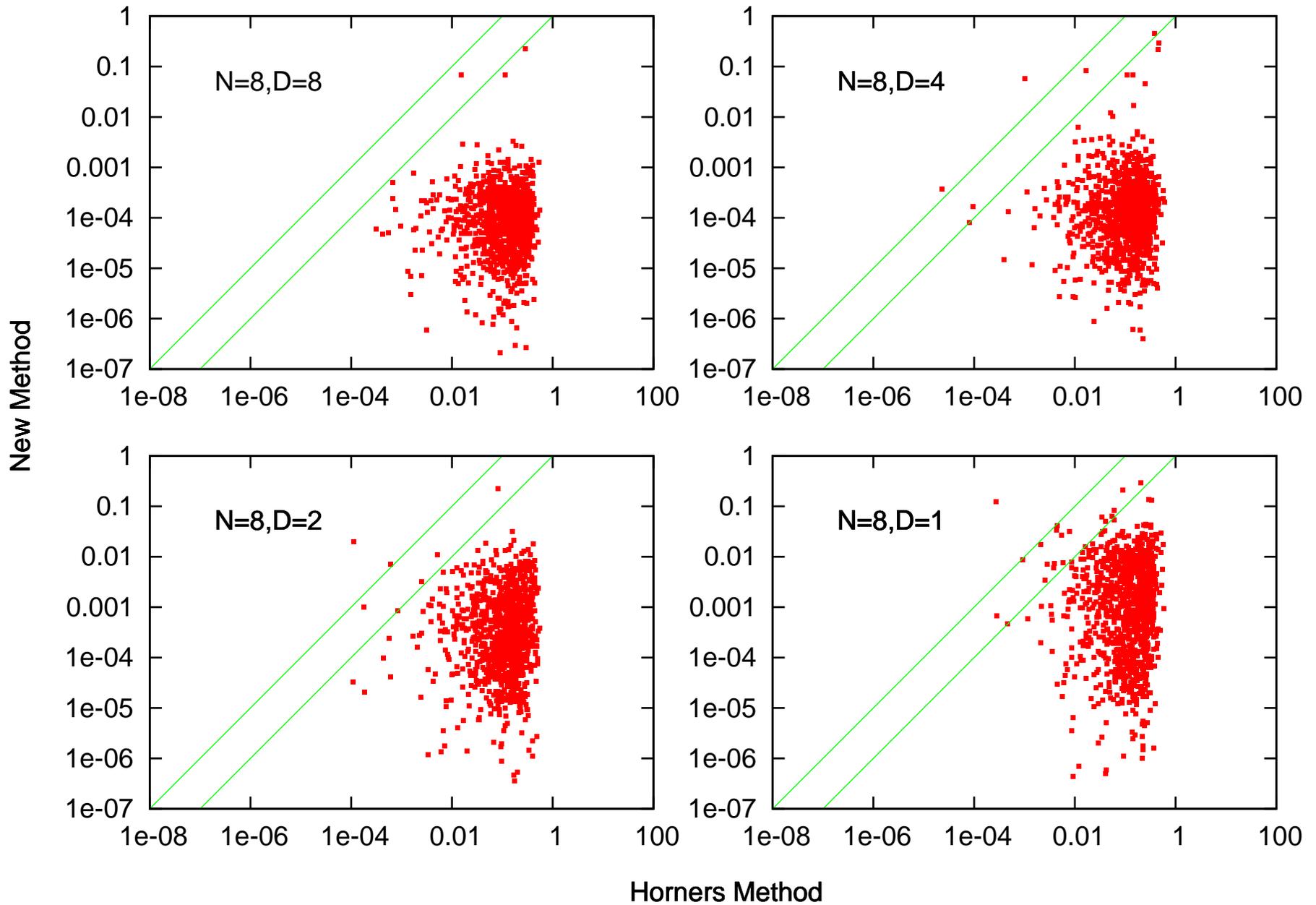

Relative Error for 8th Order Random Polynomials

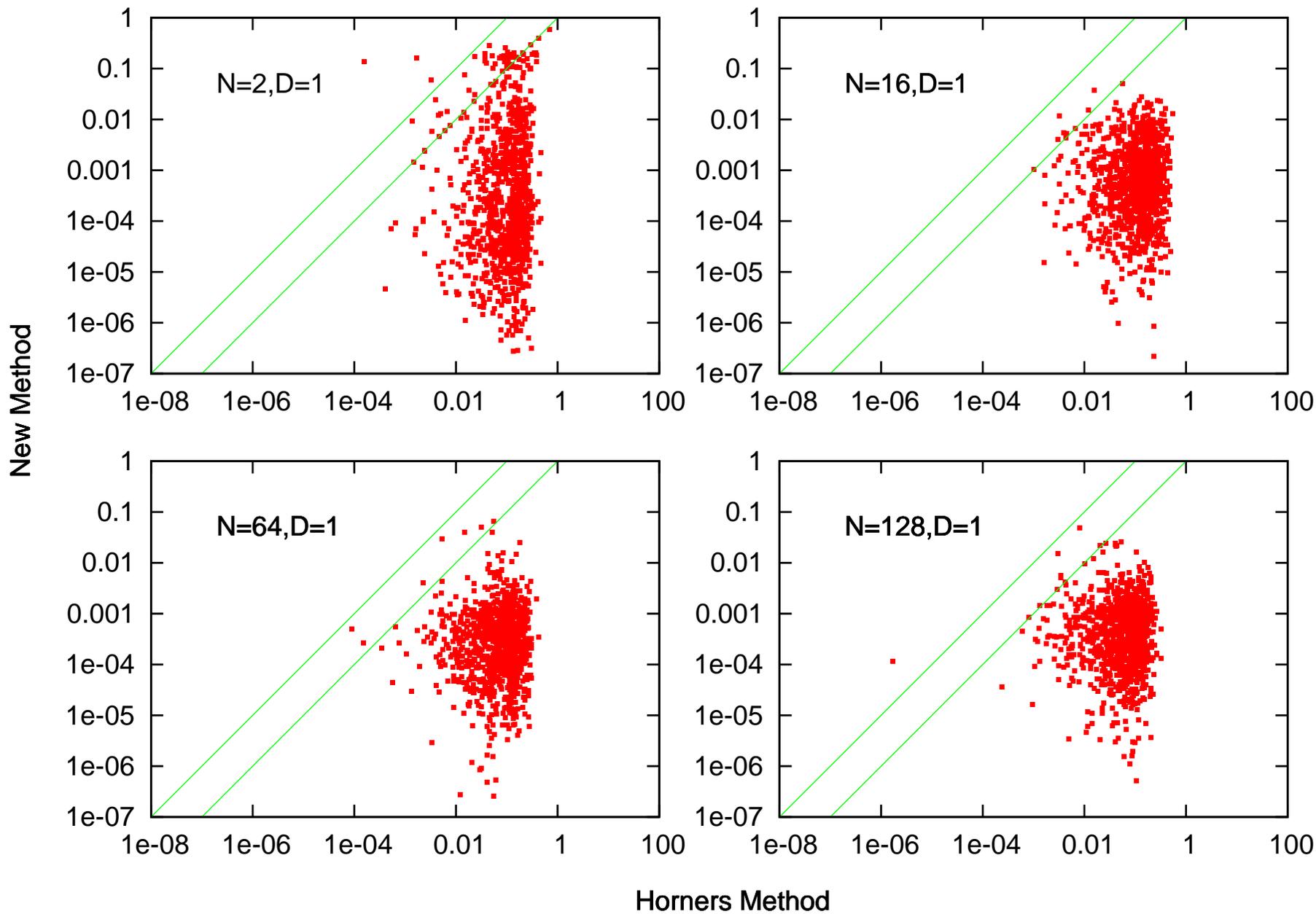